\documentclass[10pt,draft]{amsart}
\usepackage{amsmath,amssymb,amsthm}

\begin{document}
\newtheorem{lem}{Lemma}[section]
\newtheorem{prop}{Proposition}[section]
\newtheorem{cor}{Corollary}[section]
\numberwithin{equation}{section}
\newtheorem{thm}{Theorem}[section]
\theoremstyle{remark}
\newtheorem{example}{Example}[section]
\newtheorem*{ack}{Acknowledgment}
\theoremstyle{definition}
\newtheorem{definition}{Definition}[section]
\theoremstyle{remark}
\newtheorem*{notation}{Notation}
\theoremstyle{remark}
\newtheorem{remark}{Remark}[section]
\newenvironment{Abstract}
{\begin{center}\textbf{\footnotesize{Abstract}}%
\end{center} \begin{quote}\begin{footnotesize}}
{\end{footnotesize}\end{quote}\bigskip}
\newenvironment{nome}
{\begin{center}\textbf{{}}%
\end{center} \begin{quote}\end{quote}\bigskip}

\newcommand{\triple}[1]{{|\!|\!|#1|\!|\!|}}
\newcommand{\xx}{\langle x\rangle}
\newcommand{\ep}{\varepsilon}
\newcommand{\al}{\alpha}
\newcommand{\be}{\beta}
\newcommand{\de}{\partial}
\newcommand{\la}{\lambda}
\newcommand{\La}{\Lambda}
\newcommand{\ga}{\gamma}
\newcommand{\del}{\delta}
\newcommand{\Del}{\Delta}
\newcommand{\sig}{\sigma}
\newcommand{\ome}{\omega}
\newcommand{\Ome}{\Omega}
\newcommand{\C}{{\mathbb C}}
\newcommand{\N}{{\mathbb N}}
\newcommand{\Z}{{\mathbb Z}}
\newcommand{\R}{{\mathbb R}}
\newcommand{\Rn}{{\mathbb R}^{n}}
\newcommand{\Rnu}{{\mathbb R}^{n+1}_{+}}
\newcommand{\Cn}{{\mathbb C}^{n}}
\newcommand{\spt}{\,\mathrm{supp}\,}
\newcommand{\Lin}{\mathcal{L}}
\newcommand{\SSS}{\mathcal{S}}
\newcommand{\F}{\mathcal{F}}
\newcommand{\eei}{\langle\eta\rangle}
\newcommand{\xei}{\langle\xi-\eta\rangle}
\newcommand{\yy}{\langle y\rangle}
\newcommand{\dint}{\int\!\!\int}
\newcommand{\hatp}{\widehat\psi}
\renewcommand{\Re}{\;\mathrm{Re}\;}
\renewcommand{\Im}{\;\mathrm{Im}\;}

\title[On the Local Smoothing for 
the Schr\"odinger Equation]%
{On the Local Smoothing for 
the Schr\"odinger Equation}

\author{Luis Vega}

\address{Luis Vega\\
Universidad del Pais Vasco, Apdo. 64\\
48080 Bilbao, Spain}

\email{mtpvegol@lg.ehu.es}


\author{Nicola Visciglia}

\address{Nicola Visciglia\\
Dipartimento di Matematica Universit\`a di Pisa\\
Largo B. Pontecorvo 5, 56100 Pisa, Italy}

\email{viscigli@mail.dm.unipi.it}

\thanks{\noindent This research was supported by HYKE (HPRN-CT-2002-00282).
The first author was supported also by a MAC grant (MTM 2004-03029) 
and the second one by an INDAM (Istituto Nazionale di Alta Matematica) fellowship}

\date{}


\begin{abstract} We prove a family of identities 
that involve the solution $u$ to the following Cauchy problem: 
\begin{equation*}
{\bf i} \partial_t u + \Delta u=0, u(0)=f(x), (t, x)\in 
{\mathbf R}_t\times {\mathbf R}^n_x,
\end{equation*}
and the $\dot H^\frac 12({\mathbf R}^n)$-norm
of the initial datum $f$.
As a consequence of these identities we
shall deduce a lower bound for the local smoothing estimate 
proved in \cite{SC}, \cite{Sj} and \cite{V}  
and a uniqueness
criterion for the solutions to the Schr\"odinger equation.
\end{abstract}

\maketitle
\vspace{0.2cm}

This paper is devoted to the study of the following Cauchy problem:
\begin{equation}\label{schr}
{\bf i}\partial_t u + \Delta u=0, u(0)=f(x), (t, x)\in \mathbf R_t \times \mathbf R^n_x, n\geq 1.
\end{equation}

It is well-known that the solution to
\eqref{schr} satisfies the 
local smoothing estimate (see \cite{SC}, \cite{Sj} and \cite{V}):
\begin{equation}\label{smoothing}
\sup_{R\in (0, \infty)} \frac{1}{R} \int_{-\infty}^\infty \int_{B_R}
|\nabla u |^2dx dt \leq C \|f\|_{\dot H^\frac 12({\mathbf R}^n)}^2 \hbox{  }
\forall f\in \dot H^\frac 12({\mathbf R}^n),
\end{equation}
where 
$B_R$ denotes the ball in ${\mathbf R}^n$ centered in the origin of radius $R$,
$\nabla$ denotes the gradient with respect to the space variables
and $\dot H^\frac 12(\mathbf R^n)$ is the usual
homogeneous Sobolev space.

\vspace{0.1cm}

Let us recall that the estimate \eqref{smoothing}
has played a crucial role in 
the study of the
nonlinear Schr\"odinger equation 
with nonlinearities involving derivatives
(see \cite{KPV}).

\vspace{0.1cm}

Some questions can be raised in connection
with 
the local smoothing stated above. 
It is natural to ask whether
the l.h.s. in \eqref{smoothing} can be bounded from below
in the following way:
\begin{equation}\label{reversesmoothing}
\sup_{R\in (0, \infty)} \frac{1}{R} \int_{-\infty}^\infty \int_{B_R}
|\nabla u |^2dx dt \geq c \|f\|_{\dot H^\frac 12({\mathbf R}^n)}^2,
\end{equation}
where as usual $u(t, x)$ solves \eqref{schr} with initial datum $f$
and $c>0$ is a suitable constant.
Notice that an estimate of this type 
implies that \eqref{smoothing}
is an equivalence more than an inequality.

\vspace{0.1cm}

Another natural question connected with \eqref{smoothing} concerns the behaviour at infinity of
the following function:
$$F_f: (0, \infty)\ni R \rightarrow \frac 1R \int_{-\infty}^\infty \int_{B_R} 
|\nabla u|^2 dxdt\in (0, \infty),$$
where as usual $u(t, x)$ is the unique solution to \eqref{schr}.
In fact, in the best of our knowledge, it is not known 
whether the following implication is true:

\vspace{0.1cm}

\begin{equation}\label{UCP}\liminf_{R\rightarrow \infty} F_f(R)=0
\Longrightarrow f=0. 
\end{equation}

\vspace{0.1cm}

Notice that a positive answer
to \eqref{UCP} gives a uniqueness criterion
for the solutions to
the Schr\"odinger equation.

As a by product of the 
results of this paper, we can deduce that
\eqref{reversesmoothing} and \eqref{UCP}
are true.

\vspace{0.1cm}

In order to state our basic result
we have to fix some notations.

\vspace{0.1cm}

\noindent{\bf Notation.} For any $s\in {\mathbf R}$ and for any $n\in {\mathbf N}$,
the spaces $\dot H^s({\mathbf R}^n)$ shall denote the homogeneous
Sobolev spaces of order $s$,
whose norm is defined as follows:
$$\|f\|_{\dot H^s}^2=\int_{{\mathbf R}^n} |\hat f(\xi )|^2|\xi|^{2s}
d\xi,$$
where 
$$\hat f(\xi):=\int_{{\mathbf R}^n}e^{-2\pi 
{\bf i} x\xi} f(x) dx.$$

In 
the case $s=0$ we shall also use the
notation $\dot H^0({\mathbf R}^n)=L^2({\mathbf R}^n)$.

We shall denote by ${\mathcal S}({\mathbf R}^n)$
the Schwartz functional space.

If $f\in C^\infty({\mathbf R}^n)$, then  
$\partial _r f$ and $\nabla_\tau f$ 
denote respectively the radial derivative of $f$ 
and the tangential part of the full gradient $\nabla f$.

If $z=x + {\bf i} y\in {\mathbf C}$ is a complex number,
then   
${\mathcal Re} \hbox{  } z$ and ${\mathcal Im} \hbox{ }z$
shall denote its real and imaginary part.
We shall denote by $\bar z$ its complex conjugate number,
i.e. $\bar z=x- {\bf i} y$.

For any $R>0$ we shall denote by $B_R$ the ball of 
${\mathbf R}^n$ centered in the origin of radius $R$.

\vspace{0.3cm}

We can now state the basic theorem of this paper.

\begin{thm}\label{main}
Let $\psi\in C^\infty({\mathbf R}^n)$ be a real-valued
and radially symmetric function such that:
\begin{enumerate}
\item the following estimates hold:
$$ |\partial_r \psi(|x|)|, |\partial^2_r \psi(|x|)|, |\partial^3_r \psi(|x|)|
\leq C \left |P(|x|)\right | \hbox{  } \forall |x|\in {\mathbf R}^+$$
where $C>0$ and $P(|x|)$ is
a polynomial;
\item
the following limit exists:
\begin{equation*}
\lim_{|x|\rightarrow \infty} \partial_{r} \psi(|x|)
:=\psi'(\infty)\in (-\infty, \infty).
\end{equation*}\\
\end{enumerate}
Then we have the following identity:
\begin{equation}\label{global}
\lim_{T\rightarrow \infty}
\int_{-T}^T\int_{{\mathbf R}^n} \left[\nabla \bar u(t, x) D^2 \psi(x) \nabla u(t, x)
-\frac 14  |u(t, x)|^2\Delta^2 \psi(x)  \right ]dt dx
\end{equation}
\begin{equation*}=2\pi \psi'(\infty) \|f\|_{\dot H^\frac 12
({\mathbf R}^n)}^2 \hbox{  }
\forall f\in {\mathcal S}({\mathbf R}^n),
\end{equation*}
where
$u(t, x)$ is the unique solution to \eqref{schr}
with initial datum $f$, $D^2 \psi$ is the hessian matrix
$\left( \frac{\partial^2 \psi}{\partial x_i \partial x_j} \right)_{i,j=1,..,n}$
and $\Delta^2$ denotes the bilaplacian operator.
\end{thm}

\begin{remark}
Let us point out that already the existence of the limit
in the l.h.s. in \eqref{global}
is not a trivial fact, thus its existence must be considered
as a part of the statement.
\end{remark}
\begin{remark}
Let us recall that in \cite{BRV} the authors were able to show an inequality
between the l.h.s. and the r.h.s. of \eqref{global}.
Then, the main point in \eqref{global} is that 
it represents an identity
and not only an inequality.
\end{remark}

\begin{remark}
The identity \eqref{global}
can be exploited in many ways.
One possibility is to choose in \eqref{global}
a function $\psi$ such that: it satisfies
all the assumptions of the theorem \ref{main},
it is convex and  $\Delta^2 \psi(x)\leq 0$ 
for any $x\in {\mathbf R}^n$.
This was the main strategy used
in \cite{BRV} in order
to prove the local smoothing estimate
\eqref{smoothing} 
in dimension $n\geq 3$
and in presence of a potential type perturbation.
Notice that in this way 
it is possible to give a proof of \eqref{smoothing}, 
at least in dimension $n\geq 3$, which does not involve the Fourier 
transform, that was a basic tool in \cite{SC}, \cite{Sj}
and \cite{V}.
However in the sequel we shall exploit
\eqref{global} in a different direction
by taking the advantage of the fact that
it is an identity.
\end{remark}

\begin{remark}\label{remLP}
Let us notice that the identities proved in \cite{LP}
in dimension $n\geq 3$, follow from the general identity
\eqref{global} by choosing $\psi(|x|)=|x|$.
It is clear that this choice for the function $\psi(|x|)$
is not a-priori allowed since $x\rightarrow |x|$ is not a $C^\infty$ function.
However to overcome this difficulty we can choose in \eqref{global} the functions
$\psi(|x|)=\psi_{\epsilon}(|x|)=\sqrt {\epsilon^2 + |x|^2}$
and to take the limit for $\epsilon\rightarrow 0$ in the corresponding identities.
\end{remark}

As consequence of theorem \ref{main} we get
the following
result in the spirit of those given in \cite{AH}.

\begin{cor}\label{Agmon}
Assume that $n\geq 1$ and
$u(t, x)$ solves \eqref{schr} with initial datum $f$,
then the following identity holds:
\begin{equation}\label{radialident}
\lim_{R\rightarrow \infty}\frac 1R \int_{-\infty}^\infty \int_{B_R} 
|\partial_r u|^2 dxdt= 2\pi \|f\|_{\dot H^\frac 12({\mathbf R}^n)}^2 \hbox{  }
\forall f\in \dot H^\frac 12({\mathbf R}^n),
\end{equation}
and in particular
\begin{equation}\label{reversesmoothprecise}
\sup_{R\in (0, \infty)}\frac{1}{R}\int_{-\infty}^\infty \int_{B_R}|\nabla u|^2 dxdt
\geq 2\pi \|f\|_{\dot H^\frac 12({\mathbf R}^n)}^2 \hbox{ } 
\forall f\in \dot H^\frac 12({\mathbf R}^n).
\end{equation}
If moreover we assume that $u(t,x)$ satisfies the following condition:
\begin{equation}
\label{scip}\liminf_{R\rightarrow \infty} \frac 1R 
\int_{-\infty}^\infty \int_{B_R} |\nabla u|^2 dxdt=0,
\end{equation}
then $u(t, x)=0$.
\end{cor}

\begin{remark}
Notice that
the existence of the limit in \eqref{radialident}
is not a trivial fact and it must be considered
as a part of the statement.
\end{remark}

The rest of the paper is organized as follows:
section \ref{secmain}
is devoted to the proof of theorem \ref{main} while in
section \ref{agmonsec} we shall prove corollary \ref{Agmon}.

\vspace{0.3cm}

{\em {\bf Acknowledgment.}
The first author is grateful to Scuola Normale Superiore
and Centro de Giorgi in Pisa and is partially  supported by the grant MTM2005-08430 of MEC 
(Spain) and FEDER.
The second one thanks the Department of Mathematics of
the Universidad del Pais Vasco in Bilbao and INDAM 
(Istituto Nazionale di Alta Matematica). Both authors are partially supported by the European Project HPRN-CT-2002-00282-HYKE.}

\section{Proof of Theorem \ref{main}}\label{secmain}

Let us start this section with the following

\begin{lem}\label{asintotico}
If $u(t,x)$ is the unique solution to \eqref{schr} 
where $f\in {\mathcal S}({\mathbf R}^n)$, and $\psi(|x|)$ is as in theorem
\ref{main}, then:
\begin{equation}\label{radiation}
\lim_{t\rightarrow  \pm \infty}
{\mathcal Im} \int_{{\mathbf R}^n} \bar u(t, x) \nabla 
\psi(x) \nabla u(t, x) dx=\pm 2\pi \psi'(\infty) 
\|f\|_{\dot H^\frac 12({\mathbf R}^n)}^2. 
\end{equation}
\end{lem}

\noindent {\bf Proof.}
In the proof we shall need the following asymptotic formula
for the solution $u(t, x)$
to \eqref{schr} with initial datum $f\in {\mathcal S}({\mathbf R}^n)$:
\begin{equation}\label{asymptotique+}
\lim_{t\rightarrow \pm \infty} \left \|u(t,x) - 
e^{\mp {\bf i}n\pi/4}\frac{e^{\pm {\bf i} \frac{|x|^2}
{4t}}}{(4\pi t)^{n/2}} 
\hat f \left( \pm \frac x{4\pi |t|} \right)
\right\|_{L^2({\mathbf R^n})}=0,
\end{equation}
see \cite{RS} for a proof.

We first compute the limit in \eqref{radiation} as $t\rightarrow +\infty$.
By using \eqref{asymptotique+} we get
$$\lim_{t\rightarrow \infty} \|u(t,x)-v(t,x)
\|_{L^2({\mathbf R}^n)}=0,$$
where $v(t,x):=e^{-{\bf i}n\pi/4}\frac{e^{{\bf i} \frac{|x|^2}{4t}}}
{(4\pi t)^{n/2}} 
\hat f \left( \frac x{4\pi t} \right)$.

On the other hand if $u$ satisfies \eqref{schr}, then 
its partial derivatives $\partial_j u$
are still solutions of \eqref{schr} 
with Cauchy data $\partial_j f$,
for any $j=1,...,n$.

We can then apply again \eqref{asymptotique+} 
in order to deduce the following fact:
$$\lim_{t\rightarrow \infty} \|\partial_j u(t, x)- 
w_j(t, x)\|_{L^2(\mathbf R^n)}=0 \hbox{  } \forall j=1,...,n,$$
where
$$w_j (t, x):= e^{-{\bf i}n\pi/4}\frac{e^{{\bf i} 
\frac{|x|^2}{4t}}}{(4\pi t)^{n/2}} 
{\bf i} \frac{x_j}{2 t}\hat f\left( \frac x{4
\pi t} \right).$$

Notice that we have used the identity $\widehat {\partial_j 
f}(\xi)={2\pi \bf i}\xi_j \hat f(\xi)$.
We can now easily deduce that
\begin{equation*}\lim_{t\rightarrow \infty}
\int_{{\mathbf R}^n} [\bar u(t, x) \partial_j  u(t, x) 
- \bar v(t, .)  w_j(t, .)] \phi_j(x)dx=0
\end{equation*}
\begin{equation*} 
\forall \phi_j\in L^\infty({\mathbf R}^n), j=1,...,n.
\end{equation*}

In particular if we choose $\phi_j(x)=\partial_{r} \psi( |x|) \frac{x_j}{|x|}$, then we have:
$$\lim_{t\rightarrow \infty} {\mathcal Im}\int
_{{\mathbf R}^n} \bar u(t, x) \nabla \psi(x) \nabla u(t, x) dx
$$$$=\lim_{t\rightarrow \infty}
{\mathcal Im} \sum_{j=1}^n \int_{{\mathbf R}^n} \bar u(t,x) \partial_j u(t, x) 
\frac{x_j}{|x|}  \partial_{r} \psi( |x|) dx$$$$
=\lim_{t\rightarrow \infty}
{\mathcal Im} \sum_{j=1}^n \int_{{\mathbf R}^n} \bar v(t, x) 
w_j(t, x) \frac{x_j}{|x|}  \partial_{r} \psi( |x|) dx  
$$$$=
\lim_{t\rightarrow \infty}
{\mathcal Im}  \frac{{\bf i}}{(4\pi t)^n} 
\int_{{\mathbf R}^n} \frac{|x|}{2t} \left
|\hat f \left ( \frac x{4\pi t}\right )
\right |^2\partial_{r} \psi( |x|)dx.
$$

The previous chain of identities, combined with the change of variable formula,
imply:
$$\lim_{t\rightarrow \infty} {\mathcal Im}
\int_{{\mathbf R}^n} \bar u(t,x) \nabla \psi(x) \nabla u(t,x) dx$$
$$=
2\pi \lim_{t\rightarrow \infty}\int_{{\mathbf R}^n} |y|
| \hat f (y)|^2 \partial_{r} \psi( 4\pi t |y|) dy=
2\pi \psi'(\infty) \int_{{\mathbf R}^n} |y|
| \hat f (y) |^2 dy.
$$

\vspace{0.1cm}

The limit as $t\rightarrow -\infty$ in \eqref{radiation}
can be computed in a similar way by exploiting
\eqref{asymptotique+} in the case
$t\rightarrow -\infty$.

\hfill$\Box$

\noindent {\bf Proof of theorem \ref{main}.}
Following \cite{BRV} we
multiply \eqref{schr} by
\begin{equation*}\nabla \psi \nabla \bar u
+ \frac 12 \Delta \psi \bar u
\end{equation*} 
and we integrate by parts.

Let us start by writing the following identities:

$${\bf i} \partial_t u\left (\frac 12 \Delta \psi \bar u + \nabla \psi \nabla \bar u\right )$$
$$=\frac {\bf i}2  \left [div (\nabla \psi \bar u \partial_t u)
-\bar u \nabla \psi \nabla \partial_t u - \partial_tu 
\nabla \psi \nabla \bar u \right ]
$$
$$+ {\bf i}\partial_t (u \nabla \psi \nabla \bar u) - {\bf i} u \nabla \psi \nabla \partial_t \bar u$$
$$= \frac {\bf i}2 div (\nabla \psi \bar u \partial_t u) +{\bf i} 
\partial_t (u \nabla \psi \nabla \bar u)
$$
$$- {\bf i} \left ( \frac 12 \bar u \nabla \psi \nabla \partial_t u + \frac 12 \partial_tu 
\nabla \psi \nabla \bar u + u \nabla \psi \nabla \partial_t \bar u\right ).$$

Taking the real part in the previous identity we get:
$${\mathcal Re} \hbox{  } {\bf i}\partial_t u
\left (\frac 12 \Delta \psi \bar u + \nabla \psi \nabla \bar u\right )$$
$$= {\mathcal Re}\hbox{  }
{\bf i} \left [\frac 12 div (\nabla \psi \bar u \partial_t u) + 
\partial_t (u \nabla \psi \nabla \bar u)\right ]$$
$$+ {\mathcal Im} \hbox{  }
\left (\frac 12 \bar u \nabla \psi \nabla \partial_t u + \frac 12 \partial_tu 
\nabla \psi \nabla \bar u + u \nabla \psi \nabla \partial_t \bar u\right )$$
$$= {\mathcal Re}\hbox{  }
{\bf i} \left [\frac 12 div (\nabla \psi \bar u \partial_t u) + 
\partial_t (u \nabla \psi \nabla \bar u)\right ]
$$ 
$$
+ \frac 12 {\mathcal Im}\hbox{  }(\partial_tu 
\nabla \psi \nabla \bar u + u \nabla \psi \nabla \partial_t \bar u)
$$
$$= {\mathcal Re}\hbox{  }
{\bf i} \left [\frac 12 div (\nabla \psi \bar u \partial_t u) + 
\partial_t (u \nabla \psi \nabla \bar u)\right ]+ 
\frac 12 {\mathcal Im}\hbox{  } \partial_t (u \nabla \psi \nabla \bar u)
$$

If we integrate this identity on the strip $(-T, T)\times {\mathbf R}^n$ 
and we use the divergence theorem together with the assumptions done
on the growth of the derivatives of $\psi$, then we get:
\begin{equation}\label{primaparte}{\mathcal Re}\hbox{  } \int_{-T}^T\int_{{\mathbf R}^n}
{\bf i}\partial_t u\left (\frac 12 \Delta \psi \bar u + \nabla \psi \nabla \bar u\right )dxdt
\end{equation}
$$=-\frac 12 {\mathcal Im} 
\hbox{ } \sum_\pm \pm \int_{{\mathbf R}^n} \nabla \psi \nabla \bar u (\pm T, .)u(\pm T, .) dx.$$

On the other hand we have:
$${\mathcal Re}
\left [\Delta u \left (\nabla \psi \nabla \bar u
+ \frac 12 \Delta \psi \bar u\right )\right ]
$$
$$={\mathcal Re} \left [ div \left (\nabla u (\nabla \psi \nabla \bar u)
\right) - \nabla u \nabla( \nabla \psi \nabla \bar u)\right .$$
$$\left . + \frac 12 div \left (\nabla u (\Delta \psi(x) \bar u
)\right) - \frac 12 \nabla u \nabla (\bar u \Delta \psi )\right ]
$$
$$={\mathcal Re} \left [ div \left (\nabla u (\nabla \psi \nabla \bar u)
\right) + \frac 12 div \left (\nabla u (\Delta \psi \bar u)
\right) \right .$$

$$\left . -\nabla u D^2\bar u \nabla \psi  - \nabla u D^2\psi \nabla \bar u
-\frac 12 |\nabla u|^2 \Delta \psi - \frac 12 \bar u \left (\nabla u \nabla (\Delta \psi)
\right) \right ]$$

$$={\mathcal Re} \left [ div \left (\nabla u (\nabla \psi \nabla \bar u
)\right) + \frac 12 div \left (\nabla u (\Delta \psi \bar u
)\right) \right .$$

$$ 
-\frac 12 \nabla (|\nabla u|^2) \nabla \psi
-\left . \nabla u D^2\psi \nabla \bar u
-\frac 12 |\nabla u|^2 \Delta \psi -\frac 14 \nabla (|u|^2) \nabla (\Delta \psi) \right ]$$ 

$$={\mathcal Re} \left [ div \left (\nabla u (\nabla \psi \nabla \bar u)
\right) + \frac 12 div \left (\nabla u (\Delta \psi(x) \bar u)
\right) \right .$$
$$-\frac 12 div (|\nabla u|^2 \nabla \psi) 
+\frac 12 |\nabla u|^2 \Delta \psi
-\left . \nabla u D^2\psi \nabla \bar u\right. $$$$\left.
-\frac 12 |\nabla u|^2 \Delta \psi -\frac 14 div\left((|u|^2) \nabla (\Delta \psi)\right)
+ \frac 14 |u|^2 \Delta^2 \psi \right ].$$

If we integrate this identity on the strip $(-T, T)\times {\mathbf R}^n$
and we use the divergence theorem as above, then we get
\begin{equation}\label{secondaparte}
{\mathcal Re}\int_{-T}^T\int_{{\mathbf R}^n}
\Delta u \left (\nabla \psi \nabla \bar u
+ \frac 12 \Delta \psi \bar u \right ) dxdt\end{equation}
$$=
\int_{-T}^T \int_{{\mathbf R}^n} \left (- \nabla u D^2 \psi \nabla \bar u
+ \frac 14 |u|^2 \Delta^2 \psi \right )
dxdt.$$

As a consequence of the identities \eqref{primaparte} and \eqref{secondaparte},
we can deduce the following one:
\begin{equation}\label{BarcRuiVeg}
\int_{-T}^T\int_{{\mathbf R}^n} \left(\nabla \bar u D^2 \psi \nabla u
-\frac 14  |u |^2\Delta^2 \psi  \right )dt dx\end{equation}
\begin{equation*} = - \frac 12
{\mathcal Im}
\sum_\pm \pm \int_{{\mathbf R}^n} u(\pm T, .) \nabla \psi \nabla \bar u(\pm T, .) dx. 
\end{equation*}

By taking the limit as $T\rightarrow \infty$ in \eqref{BarcRuiVeg} 
and by using lemma 
\ref{asintotico}, we can deduce the desired result.

\hfill$\Box$

\section{Applications}\label{agmonsec}

This section is devoted to the proof of corollary
\ref{Agmon}.

In order to do to that
we shall need the following 
lemma.

\begin{lem}\label{noma}
Assume that $f\in {\mathcal S}({\mathbf R}^n)$
is such that: 
\begin{equation}\label{dense}\hat f(\xi)=0 \hbox{  } 
\forall \xi=(\xi_1,...,\xi_n)\in \mathbf R^n_\xi
\hbox{ s.t. } |\xi_1|<\epsilon_0=\epsilon_0(f) \hbox{ where }\epsilon_0>0.
\end{equation} 
If $u(t,x)$ is the corresponding solution
to \eqref{schr},
then:
\begin{equation}\label{zero}
\lim_{R\rightarrow \infty}  \int_{-\infty}^\infty\int_{{\mathbf R}^n}
\frac{|\nabla_\tau u (t, x)|^2}{|x|} |\partial_r \phi_R(|x|)| 
dxdt =0,
\end{equation}
\begin{equation}\label{uno}
\lim_{R\rightarrow \infty}  \int_{-\infty}^\infty\int_{{\mathbf R}^n}
|u(t, x)|^2 |\Delta^2 \phi_R(|x|)| 
dxdt =0,\end{equation}
where $\phi(|x|)\in C^\infty({\mathbf R}^n)$ is any  radially symmetric function
such that: $$\partial_{r}\phi(|x|)\leq C, |\Delta^2 \phi (|x|)|\leq \frac{C}
{(1+|x|)^3} \hbox{  } \forall x\in {\mathbf R}^n,$$
and $\phi_R(|x|):=R \phi\left (\frac{|x|}R\right )$.
\end{lem}

\noindent{\bf Proof.}
Notice that \eqref{zero} follows by combining 
the following implication:
$$u \hbox{ solves } \hbox{ \eqref{schr} } \hbox{ with } f\in \mathcal{S}({\mathbf R}^n)
\Rightarrow \frac{|\nabla_\tau u|^2}{|x|}
\in L^1({\mathbf R}_t \times {\mathbf R}_x^n),$$
whose proof can be found in \cite{LP} (see also remark \ref{remLP}),
with the following trivial fact:
$$\lim_{R\rightarrow \infty} \partial_r \phi_R(|x|)=\partial_r \phi(0)=0
\hbox{  } \forall x\in {\mathbf R}^n,$$
where we have used the radiality of $\phi(|x|)$.

\vspace{0.1cm}

Next we shall show \eqref{uno}.
Let us introduce the unique function $g\in {\mathcal S}({\mathbf R}^n)$ such that:
$$\hat g(\xi)=\frac{\hat f(\xi)}{2\pi {\bf i}\xi_1} \hbox{  } \forall \xi \in {\mathbf R}^n_\xi$$
and let us consider the unique solution $v(t,x)$ to
\eqref{schr} with inital datum given by $g(x)$.

It is easy to check that $\partial_{x_1}v(t,x)=u(t,x)$.
We can then apply \eqref{smoothing} to the solution $v(t, x)$ in order to get:
\begin{equation}\label{colper}
\frac 1R \int_{-\infty}^\infty\int_{B_R} |u(t,x)|^2 dxdt
=\frac 1R \int_{-\infty}^\infty\int_{B_R} |\partial_{x_1}v(t,x)|^2 dxdt
\end{equation}
$$
\leq C  \|g \|_{\dot H^\frac 12({\mathbf R}^n)}^2
<\infty \hbox{  } \forall R>0.$$

Notice that due to the assumption done on $\phi(|x|)$ we get:
$$\int_{-\infty}^\infty \int_{{\mathbf R}^n} |u(t,x)|^2 
|\Delta^2 \phi_R(|x|)|dx dt = \frac{1}{R^3} \int_{-\infty}^\infty \int_{{\mathbf R}^n} |u(t,x)|^2 
\left |\Delta^2 \phi\left (\frac{|x|}{R}\right )\right |dx dt$$
$$\leq C \left ( \frac 1{R^3} \int_{-\infty}^\infty \int_{|x|<1} |u(t,x)|^2dxdt +\sum_{j=1}^\infty
\int_{-\infty}^\infty \int_{2^j< |x|<2^{j+1}} \frac{|u(t,x)|^2}{(2^j + R)^3} dxdt \right)
$$
$$\leq C \left ( \frac 1{R^3}\int_{-\infty}^\infty \int_{|x|<1} |u(t,x)|^2dxdt +
\sum_{j=1}^\infty \int_{-\infty}^\infty \int_{|x|<2^{j+1}} \frac{|u(t,x)|^2}{(2^j + R)^3} dxdt \right),
$$
that due to \eqref{colper} implies:
\begin{equation}\label{tart} \int_{-\infty}^\infty \int_{{\mathbf R}^n} |u(t,x)|^2 
|\Delta^2 \phi_R(|x|)|dx dt\end{equation}
$$
\leq C \left (\frac 1{R^3}+  \sum_{j=0}^\infty \frac{2^{j+1}}{(2^j+R)^3}\right ).$$

On the other hand it is easy to show that
$$\lim_{R\rightarrow \infty} \frac 1{R^3}+  \sum_{j=0}^\infty \frac{2^{j+1}}{(2^j+R)^3}
=0,$$
that in conjunction with \eqref{tart}
implies
\eqref{uno}.

\hfill$\Box$

\vspace{0.1cm}

\begin{remark}
Let us notice that \eqref{uno}
follows, at least in dimension $n\geq 4$,
from the following inequality:
\begin{equation}\label{LPdo}\int_{-\infty}^\infty \int_{{\mathbf R}^n}\frac {|u|^2}{|x|^3} dxdt
\leq C\|f\|_{\dot H^\frac 12({\mathbf R}^n)} \hbox{   } \forall n\geq 4
\end{equation}
whose proof can be found in \cite{LP}. 
\end{remark}

\vspace{0.1cm}

\noindent{\bf Proof of corollary \ref{Agmon}.}
\noindent 
First of all
we show that
\eqref{radialident} implies \eqref{reversesmoothprecise}:
$$\sup_{R\in (0,\infty)}
\frac 1 R 
\int_{-\infty}^\infty \int_{B_R}
|\nabla u|^2 dxdt \geq
\sup_{R\in (0,\infty)}
\frac 1 R 
\int_{-\infty}^\infty \int_{B_R}
|\partial_r u|^2 dxdt $$$$\geq
\lim_{R\rightarrow \infty} \frac 1 R 
\int_{-\infty}^\infty \int_{B_R}
|\partial_r u|^2 dxdt =2\pi \|f\|_{\dot H^\frac12({\mathbf R}^n)}^2.
$$ 

On the other hand if we assume 
\eqref{scip}, then by
\eqref{radialident} we get:
$$
0=\liminf_{R\rightarrow \infty}
\frac 1 R 
\int_{-\infty}^\infty \int_{B_R}
|\nabla u|^2 dxdt \geq \liminf_{R\rightarrow \infty}
\frac 1 R 
\int_{-\infty}^\infty \int_{B_R}
|\partial_r u|^2 dxdt $$$$ 
= \lim_{R\rightarrow \infty}
 \frac 1 R 
\int_{-\infty}^\infty \int_{B_R}
|\partial_r u|^2 dxdt =2 \pi \|f\|_{\dot H^\frac12({\mathbf R}^n)}^2, 
$$
then $f=0$ and in particular $u(t, x)=0$.

\vspace{0.1cm}

Next we shall prove \eqref{radialident} assuming that the initial datum $f$ is such that
$f\in \mathcal{S}({\mathbf R}^n)$ and moreover it satisfies condition \eqref{dense}.
It is easy to show, by combining a density argument
with \eqref{smoothing}, that this regularity assumption
done on $f$
can be removed.

\vspace{0.1cm}

For any $k\in \mathbf N$
we fix a function
$h_k(r)\in C^\infty_0(\mathbf R; [0, 1])$ such that:
$$h_k(r)=1 \hbox{  } \forall r\in {\mathbf R}
\hbox{ s.t. } |r|<1, h_k(r)=0 \hbox{  } \forall r \in {\mathbf R}
\hbox{ s.t. } |r|>\frac{k+1}{k}, 
$$
$$h_k(r)= h_k(-r) \hbox{  } \forall r \in \mathbf R.$$

Let us introduce the functions
$\psi_k(r), H_k(r) \in C^\infty( {\mathbf R})$:
$$\psi_k(r) = \int_0^r (r-s) h_k(s) ds \hbox{  } \hbox{ and } \hbox{  }
H_k(r):=\int_0^r h_k(s) ds.$$

Notice that 
\begin{equation}\label{differ}
\psi_k''(r)=h_k(r), \psi_k'(r)= H_k(r) \forall r\in {\mathbf R} \hbox{ and }
\lim_{r\rightarrow \infty} \partial_r \psi_k(r)=\int_0^\infty h_k(s)ds.
\end{equation}

Moreover an elementary computation 
shows that: 

\begin{equation}\label{bilaplac}\Delta^2 \psi_k(|x|)= 
\frac C{|x|^3} \hbox{  } \forall x\in 
{\mathbf R}^n \hbox{ s.t. } |x| \geq 2 \hbox{ and } n\geq 2,
\end{equation}
where $\Delta^2$ is the bilaplacian operator,
while in the one dimensional case, i.e. for $n=1$, we have: $$\partial_x^4 \psi_k(|x|)=0 \hbox{  } \forall x\in 
{\mathbf R} \hbox{ s.t. } |x| \geq 2.$$

Thus the functions
$\phi=\psi_k(|x|)$ satisfy the assumptions of lemma \ref{noma}
in any dimension $n\geq 1$.

\vspace{0.1cm}

In the sequel we shall need the rescaled functions
$$\psi_{k, R}(|x|):=R \psi_k\left (\frac{|x|}{R} \right ) \forall x\in {\mathbf R}^n, k\in 
{\mathbf N} \hbox{ and }R>0$$
and we shall exploit the following 
elementary identity:
$$\nabla \bar u D^2 \psi \nabla u=
\partial_r^2\psi(|x|) |\partial_r u|^2 + \frac{\partial_r \psi(|x|)}{|x|}|\nabla_\tau u|^2,$$ 
where $\psi(|x|)$ is any regular radial function
and $u$ is another regular function.

By combining this identity with \eqref{global}, where we choose $\psi(|x|)= 
\psi_{k,R}(|x|)$, and recalling \eqref{differ} we get:
\begin{equation}\label{limit}\int_{-\infty}^\infty \int_{{\mathbf R}^n}
\left [ \partial_r^2 \psi_{k,R} |\partial_r u|^2 + \frac{\partial_r \psi_{k,R}}
{|x|} | \nabla_\tau u|^2
-\frac 14
|u|^2 \Delta^2 \psi_{k,R}\right ]dx dt 
\end{equation}
\begin{equation*}
= 2\pi \left (\int_0^\infty h_{k}(s) ds
\right ) \|f \|_{\dot H^\frac 12 ({\mathbf R}^n)}^2 \forall k\in {\mathbf N}, R>0.
\end{equation*}

By using \eqref{zero} and \eqref{uno} where we make the choice 
$\phi(|x|)=\psi_k(|x|)$ we get
\begin{equation}\label{lim1}
\lim_{R\rightarrow \infty} \int_{-\infty}^\infty \int_{{\mathbf R}^n}
\left [\frac{|\nabla_\tau u|^2}{|x|}\partial_r \psi_{k,R} 
- \frac 14 \Delta^2 \psi_{k,R}  |u|^2 \right ]dt dx=0
\hbox{  } \forall k\in {\mathbf N}.
\end{equation}

We can combine now \eqref{limit} with \eqref{lim1} in order to deduce

\begin{equation}\label{crusia}\lim_{R\rightarrow \infty} \int_{-\infty}^\infty \int_{{\mathbf R}^n}
\partial_r^2\psi_{k,R} |\partial_r u|^2
dx dt =2\pi \left (\int_0^\infty h_{k}(s) ds
\right ) \|f \|_{\dot H^\frac 12 ({\mathbf R}^n)}^2 \forall k\in {\mathbf N}.
\end{equation}

On the other hand,
due to the properties of $h_{k}$, we get
$$\frac 1R \int_{-\infty}^\infty \int_{B_R}|\partial_r u|^2
dx dt \leq \int_{-\infty}^\infty \int_{{\mathbf R}^n}
\partial_r^2 \psi_{k,R}  |\partial_r u|^2  dt dx $$$$
= \frac 1R \int_{-\infty}^\infty \int_{{\mathbf R}^n}
h_k\left( \frac {|x|}R\right)|\partial_r u|^2  dt dx \leq \frac 1R \int_{-\infty}^\infty \int_{B_{\frac {k+1}{k}R}} |\partial_r u|^2 dxdt$$
that due to \eqref{crusia} implies:
\begin{equation}\label{limsup}\limsup_{R\rightarrow \infty} \frac 1R \int_{-\infty}^\infty \int_{B_R}|\partial_r u|^2 dxdt 
\leq 2\pi \left (\int_0^\infty h_k(s) ds \right ) \|f\|_{\dot H^\frac 12
({\mathbf R}^n)}^2\end{equation}$$\leq\frac {k+1}{k}
\liminf_{R\rightarrow \infty} \frac 1R \int_{-\infty}^\infty \int_{B_R}|\partial_r u|^2 dxdt
\hbox{  } \forall k \in {\mathbf N}.$$

Since $k\in \mathbf N$ is arbitrary and since 
the following identity is trivially satisfied:
$$\lim_{k\rightarrow \infty}\int_0^\infty h_k(s) ds=1,$$
we can deduce easily
\eqref{radialident} by using \eqref{limsup}.

\vspace{0.1cm}

The proof is complete.

\hfill$\Box$

\end{document}